\documentclass[12pt]{amsart}
\usepackage{amssymb,upref,enumerate}

\usepackage{amsmath,amssymb,amscd,amsthm,amsxtra}
\usepackage{latexsym}
\usepackage{graphics,epsfig}
\allowdisplaybreaks

\def\XXint#1#2#3{{\setbox0=\hbox{$#1{#2#3}{\int}$}
\vcenter{\hbox{$#2#3$}}\kern-.5\wd0}}

\newtheorem{theorem}{Theorem}
\newtheorem{proposition}{Proposition}

\newtheorem{lemma}{Lemma}

\newtheorem{corollary}{Corollary}

\newtheorem*{theorem*}{Theorem}
\theoremstyle{remark}

\def\({\left(}
\def\){\right)}
\def\be {\begin{equation}}
\def\en{\end{equation}}

\def\Cdot{{\dot{C}}}
\def\div{\text{div}~}

\newcommand{\tensor}{\otimes}

\numberwithin{equation}{section}

\begin{document}

%\dedicatory{}

\subjclass[2000]{Primary 76D05, 35A02; Secondary 35K58}

\keywords{Navier-Stokes, Lagrangian Averaging, global existence, Besov spaces}

%\thanks{}

\title[LANS equation in Besov spaces]{Local and global existence for the Lagrangian Averaged Navier-Stokes equations in Besov spaces}
\author{Nathan Pennington}
\address{Nathan Pennington, Department of Mathematics, Kansas State University, 138 Cardwell Hall,
Manhattan, KS-66506, USA.} \email{npenning@math.ksu.edu}

\date{\today}

\begin{abstract} 
We prove the existence of short time solutions to the incompressible, isotropic Lagrangian Averaged Navier-Stokes equation with low regularity initial data in Besov spaces $B^{r}_{p,q}(\mathbb{R}^n)$, $r>n/2p$.  When $p=2$ and $n\geq 3$, we obtain global solutions, provided the parameters $r,q$ and $n$ satisfy certain inequalities.  This is an improvement over known analogous Sobolev space results, which required $n=3$.

\end{abstract}

\maketitle

\bigskip

\section{Introduction}
The Lagrangian Averaged Navier-Stokes (LANS) equation is a recently derived approximation to the Navier-Stokes equation.  The equation is obtained via an averaging process applied at the Lagrangian level, resulting in a modified energy functional.  The geodesics of this energy functional satisfy the Lagrangian Averaged Euler (LAE) equation, and the LANS equation is derived from the LAE equation in an analogous fashion to the derivation of the Navier-Stokes equation from the Euler equation.  For an exhaustive treatment of this process, see \cite{Shkoller}, \cite{SK}, \cite{MRS} and \cite{MS2}. In \cite{MKSM} and \cite{CHMZ}, the authors discuss the numerical improvements that use of the LANS equation provides over more common approximation techniques of the Navier-Stokes equation.

On a region without boundary, the isotropic, incompressible form of the LANS equation is given by 
\begin{equation}\label{LANS}\aligned\partial_t u+(u\cdot\nabla)u+\div
\tau^\alpha u=-(1-\alpha^2\triangle)^{-1}\nabla p+\nu\triangle
u
\\ u=u(t,x),~~ \div u=0,~~ u(0,x)=u_0(x), 
\endaligned
\end{equation}
with the terms defined as follows.  First, $u:I\times \mathbb{R}^n\rightarrow \mathbb{R}^n$ for some time strip $I=[0,T)$ denotes the velocity of the fluid, $\alpha>0$ is a constant, $p:I\times\mathbb{R}^n\rightarrow \mathbb{R}^n$ denotes the fluid pressure, $\nu>0$ is a constant due to the viscosity of the fluid, and $u_0:\mathbb{R}^n\rightarrow\mathbb{R}^n$, with $\div u_0=0$.  Next, the differential operators $\nabla, \triangle,$ and $\div$ are spatial differential operators with their standard definitions.  The term $(v\cdot \nabla)w$, also denoted $\nabla_v w$, is the vector field with $j^{th}$ component $\sum_{i=1}^n v_i\partial_i w_j$.  The Reynolds
stress $\tau^\alpha$ is given by \begin{equation*}\tau^\alpha u=\alpha^2(1-\alpha^2\triangle)^{-1}[Def(u)\cdot
Rot(u)]
\end{equation*}
where $Rot(u)=(\nabla u-\nabla u^T)/2$ and $Def(u)=(\nabla u+\nabla u^T)/2$.  We remark that setting $\alpha=0$ in equation $(\ref{LANS})$ recovers the Navier-Stokes equation.  

There is a wide variety of local existence results for the LANS equation in various settings, including \cite{Shkoller}, \cite{MRS}, \cite{MS} and \cite{sobpaper}.  In \cite{MS}, Marsden and Shkoller proved the existence of global solutions to the LANS equation with initial data in the Sobolev space $H^{3,2}(\mathbb{R}^3)$.  In \cite{sobpaper}, this result was improved, achieving global existence for data in the space $H^{3/4,2}(\mathbb{R}^3)$.  
In this article we seek local and global solutions to the LANS equation with initial data in the Besov spaces $B^{s}_{p,q}(\mathbb{R}^n)$, with an emphasis placed on minimizing $s$.  These results are comparable to those in \cite{sobpaper} for Sobolev space initial data.  

The main distinction between the results in this paper and those in \cite{sobpaper} is that the extension of local results to global results in the Besov space setting does not follow from standard energy methods, requiring significantly more delicate Fourier analysis.  As a result of this argument, the global existence result can be applied in dimension $n\geq 3$, which improves the requirement $n=3$ in \cite{sobpaper}.

The paper is organized as follows.  We devote the rest of this section to defining solution spaces and stating our main theorems.  Theorem $\ref{main global piece}$ yields an $\emph{a~priori}$ bound for solutions to the LANS equation, and, as mentioned in the preceding paragraph, is the main result of the paper.  Its proof can be found in Section $\ref{global extension}$.  Theorems $\ref{special case 1}$ through $\ref{thm2}$ are local existence theorems, and Corollary $\ref{besov global extension}$ is the global existence result.  In Section $\ref{Function Space definitions and basic Besov space results}$, we recall some known Besov space operator estimates which will be useful in the proofs of our main results.  In Section $\ref{local continuous solutions}$ we prove Theorems $\ref{special case 1}$ and $\ref{thm1}$ and in Section $\ref{local integral solutions}$ we prove Theorems $\ref{special case 2}$ and $\ref{thm2}$.    

As mentioned above, we denote Besov spaces by $B^s_{p,q}(\mathbb{R}^n)$, with norm denoted by $\|\cdot\|_{B^s_{p,q}}=\|\cdot\|_{s,p,q}$ (a  complete definition of these spaces can be found in Section $2$).  We define the space
\begin{equation*}C^T_{a;s,p,q}=\{f\in
C((0,T):B^s_{p,q}(\mathbb{R}^n)):\|f\|_{a;s,p,q}<\infty\},
\end{equation*}
where
\begin{equation*}\|f\|_{a;s,p,q}=\sup\{t^a\|f(t)\|_{s,p,q}:t\in
(0,T)\},
\end{equation*}
$T>0$, $a\geq 0$, and $C(A:B)$ is the space of continuous functions from $A$ to $B$.  We let $\Cdot^T_{a;s,p,q}$ denote the subspace of $C^T_{a;s,p,q}$ consisting
of $f$ such that
\begin{equation*}\lim_{t\rightarrow 0^+}t^a f(t)=0
~\text{(in}~B^s_{p,q}(\mathbb{R}^n)).
\end{equation*}
Note that while the norm $\|\cdot\|_{a;s,p,q}$ lacks an explicit reference to $T$, there is an implicit $T$ dependence. We also say $u\in BC(A:B)$ if $u\in C(A:B)$ and $\sup_{a\in A}\|u(a)\|_{B}<\infty$.  Lastly, setting $\mathbb{M}((0,T):\mathbb{E})$ to be the set of measurable
functions defined on $(0,T)$ with values in the space $\mathbb{E}$, we define
\begin{equation*}L^a((0,T):B^s_{p,q}(\mathbb{R}^n))=
\{f\in \mathbb{M}((0,T):B^s_{p,q}(\mathbb{R}^n)):(\int_0^T
\|f(t)\|^a_{s,p,q}dt)^{1/a}<\infty\}.
\end{equation*}

As a last bit of notation, we set $r^+$ to be a real number arbitrarily close to, but strictly greater than, $r$.  Similarly, $r^-$ is arbitrarily close to, but strictly less than, $r$. 

Now we state our main theorems.  Our first result provides a uniform-in-time bound on the Besov norm of solutions to the LANS equation.

\begin{theorem}\label{main global piece}Let $u:I\times\mathbb{R}^n\rightarrow\mathbb{R}^n$ be a solution to $(\ref{LANS})$ on a time strip $I=[0,T)$.  Then
\begin{equation*}\|u(t)\|_{B^r_{2,q}}\leq \|u_0\|_{B^r_{2,q}}\exp\(C\int_I \|u(t)\|_{B^{1+n/2}_{2,q}} dt\),
\end{equation*}
where $r>2$.
\end{theorem}

We now state our four local existence theorems.  As we are working in low regularity spaces, these local solutions are not classical solutions.  Instead, they are solutions to an integral equation derived from the classical formulation of the PDE through Duhamel's principle.  A more explicit accounting of this can be found in Section $\ref{local continuous solutions}$.  Finally, all of these local solutions depend uniquely and continuously on the initial data.  

\begin{theorem}\label{special case 1}Let $s=3/4$ and let $u_0\in B^{s^+}_{2,q}(\mathbb{R}^3)$ be divergence-free.  Then there exists a local solution $u$ to the LANS equation $(\ref{LANS})$, where
\begin{equation}\label{lowry1}u\in BC([0,T):B^{s^+}_{2,q}(\mathbb{R}^3))\cap \dot{C}^T_{(1^+-s^+)/2;1^+,2,q},
\end{equation}
and $T$ is a non-increasing function of $\|u_0\|_{s^+,2,q}$, with $T=\infty$ if $\|u_0\|_{s^+,2,q}$ is sufficiently small.

Next, let $p>n\geq 3$ and let $s=n/2p$.  Then for any divergence-free $u_0\in B^{s^+}_{p,q}(\mathbb{R}^n)$, there exists a locally well-posed solution $u$ to the LANS equation $(\ref{LANS})$, where 
\begin{equation*}u\in BC([0,T):B^{s^+}_{p,q}(\mathbb{R}^n))\cap \dot{C}^T_{(1^+-s^+)/2;1^+,p,q},
\end{equation*}
and $T$ is a non-increasing function of $\|u_0\|_{s^+,2,q}$, with $T=\infty$ if $\|u_0\|_{s^+,2,q}$ is sufficiently small.
\end{theorem} 

Our third theorem is analogous to the second, with the space $L^a((0,T):B^{s}_{p,q}(\mathbb{R}^n))$ replacing $\dot{C}^T_{a;s,p,q}$.
\begin{theorem}\label{special case 2}Let $s=3/4$ and let $u_0\in B^{s_1}_{2,q}(\mathbb{R}^3)$ be divergence-free.  Then there exists a local solution $u$ to the LANS equation $(\ref{LANS})$, where 
\begin{equation}\label{lowry3}u\in BC([0,T):B^{s_1}_{2,q}(\mathbb{R}^n))\cap L^{2/(1^+-s^+)}((0,T):B^{1^+}_{2,q}(\mathbb{R}^n)),
\end{equation}
and $T$ is a non-increasing function of $\|u_0\|_{s^+,2,q}$, with $T=\infty$ if $\|u_0\|_{s^+,2,q}$ is sufficiently small.

Next, let $p>n\geq 3$ and $s=n/2p$.  Then for any divergence-free $u_0\in B^{s^+}_{p,q}(\mathbb{R}^n)$ there exists a local solution $u$ to the LANS equation $(\ref{LANS})$ where 
\begin{equation*}u\in BC([0,T):B^{s^+}_{p,q}(\mathbb{R}^n))\cap L^{2/(1^+-s^+)}((0,T):B^{1^+}_{p,q}(\mathbb{R}^n)),
\end{equation*}
and $T$ is a non-increasing function of $\|u_0\|_{s^+,2,q}$, with $T=\infty$ if $\|u_0\|_{s^+,2,q}$ is sufficiently small.
\end{theorem} 

We note that the use of $(\frac{n}{2p})^+$ in these two theorems is due to the requirement in Proposition $\ref{besovproduct}$ that $s>n(2/p_1-1/p)$ (with $2p\geq p_1$), which differs from the analogous Sobolev space requirement that $s\geq n(2/p_1-1/p)$. 

In the next two theorems, the minimal regularity achieved is $(n/p)^+$.  However, these two theorems are significantly more general, allowing for a much wider range of parameter values, specifically with regard to the auxiliary space.

\begin{theorem}\label{thm1}For any divergence-free $u_0\in B^r_{p,q}(\mathbb{R}^n)$ there exists a local solution to $(\ref{LANS})$ such that
\begin{equation}\label{lowry2}u\in BC([0,T):B^r_{p,q}) \cap \dot{C}^T_{a;s,\tilde{p},q},
\end{equation}
for some $T>0$, provided the parameters satisfy the list of conditions $(\ref{conditionsbesov1})$, including that $r<s<r+1$.  The minimal $r$ achieved by these methods is $r=(\frac{n}{p})^+$.  

Also, $T=T(\|u_0\|_{r,p,q})$ can be chosen to be a non-increasing function of $\|u_0\|_{r,p,q}$, with $T=\infty$ if $\|u_0\|_{r,p,q}$ is sufficiently small.
\end{theorem}

\begin{theorem}\label{thm2}For any divergence-free $u_0\in B^{r}_{p,q}(\mathbb{R}^n)$ there exists a local solution $u$ to $(\ref{LANS})$ such that
\begin{equation}\label{lowry4}u\in BC([0,T):B^r_{p,q})\cap
L^\sigma((0,T):B^s_{\tilde{p},q}),
\end{equation}
for some $T>0$, provided the parameters satisfy the list of conditions
$(\ref{besovthmpara})$, including that $r<s<r+1$.  The minimal $r$ achieved by these methods is $r=(\frac{n}{p})^+$.   Also, $T=T(\|u_0\|_{r,p,q})$ can be chosen to be a non-increasing function of $\|u_0\|_{r,p,q}$, with $T=\infty$ if $\|u_0\|_{r,p,q}$ is sufficiently small.
\end{theorem}

Combining these local existence results with Theorem $\ref{main global piece}$ (which requires Proposition $\ref{higher regularity theorem}$ or Proposition $\ref{higher regularity theorem2}$, depending on the context), we get the following corollary.

\begin{corollary}\label{besov global extension}The local solutions from Theorems $\ref{special case 1}$ and $\ref{thm1}$ given by equations $(\ref{lowry1})$ and $(\ref{lowry2})$ (with $p=\tilde{p}=2$), respectively, can be extended to global solutions.  The local solutions from Theorems $\ref{special case 2}$ and $\ref{thm2}$ given by equations $(\ref{lowry3})$ and $(\ref{lowry4})$ (with $p=\tilde{p}=2$), respectively, can also be extended to global solutions.
\end{corollary}

The proof of the corollary follows from standard extension methods.  

\section{Besov spaces}\label{Function Space definitions and basic Besov space results}

We begin by defining the Besov spaces $B^s_{p,q}(\mathbb{R}^n)$.  Let $\psi_0\in\mathcal{S}$ be an even, radial function such that $\hat{\psi_0}\geq 0$, the support of $\hat{\psi_0}$ lies in the annulus $A_0:=\{\xi\in \mathbb{R}^n:2^{-1}<|\xi|<2\}$, and $\sum_{j\in\mathcal{Z}} \hat{\psi_0}(2^{-j}\xi)=1$ for all $\xi\neq 0$.  We then define 
$\hat{\psi_j}(\xi)=2^{jn}\hat{\psi}_0(2^{-j}\xi)$, and remark that $\hat{\psi_j}$ is supported in $A_j:=\{\xi\in\mathbb{R}^n:2^{j-1}<|\xi|<2^{j+1}\}$.  We also define $\Psi$ by 
\begin{equation*}\hat{\Psi}(\xi)=1-\sum_{k=0}^\infty \hat{\psi}_k(\xi).
\end{equation*}

We define the operators $\triangle_j$ and $S_j$ by
\begin{equation*}\triangle_j f=\psi_j\ast f, \quad
S_jf=\sum_{k=-\infty}^{j}\triangle_k f,
\end{equation*}
and record some properties of these operators.  Applying the Fourier Transform and
recalling that $\hat{\psi}_j$ is supported on $2^{j-1}\leq
|\xi|\leq2^{j+1}$, it follows that   
\begin{equation}\aligned \label{besovlemma1}\triangle_j\triangle_k f= 0, \quad |j-k|\geq 2
\\ \triangle_j (S_{k-3}f\triangle_{k}g)= 0 \quad |j-k|\geq 3,
\endaligned
\end{equation}
and, if $|i-k|\leq 1$, then 
\begin{equation}\label{besovpieces67}\triangle_j(\triangle_kf\triangle_i g)=0 \quad j>k+3.
\end{equation}

For $s\in\mathbb{R}$ and $1\leq p,q\leq \infty$ we define
the space $\tilde{B}^s_{p,q}(\mathbb{R}^n)$ to be the set of distributions such that 
\begin{equation*}\|u\|_{\tilde{B}^s_{p,q}}=\(\sum_{j=0}^\infty (2^{js}\|\triangle_j
u\|_{L^p})^q\)^{1/q}<\infty,
\end{equation*}
with the usual modification when $q=\infty$.  Finally, we define the Besov spaces $B^s_{p,q}(\mathbb{R}^n)$ by the norm 
\begin{equation*}\|f\|_{B^s_{p,q}}=\|\Psi*f\|_p+\|f\|_{\tilde{B}^s_{p,q}},
\end{equation*}
for $s>0$.  For $s>0$, we define $B^{-s}_{p',q'}$ to be the dual
of the space $B^s_{p,q}$, where $p',q'$ are the Holder-conjugates to
$p,q$.

Now we turn our attention to establishing some basic Besov space estimates.  First, we let $1\leq q_1\leq q_2\leq \infty$, $\beta_1\leq \beta_2$, $s> 0$, $1\leq p_1\leq p_2\leq\infty$ and $\gamma_1=\gamma_2+n(1/p_1-1/p_2)$.  Then we have the following:
\begin{equation}\label{besov embedding}\aligned \|f\|_{B^{\beta_1}_{p,q_2}}&\leq C\|f\|_{B^{\beta_2}_{p,q_1}}
\\ \|f\|_{B^{\gamma_2}_{p_2,q}}&\leq C\|f\|_{B^{\gamma_1}_{p_1,q}}
\\ \|f\|_{L^p}&\leq \|f\|_{B^s_{p,q}}.
\endaligned
\end{equation}

We also have the following fractional Leibniz rule estimate (see, for instance, Lemma $2.2$ in \cite{chae}).
\begin{proposition}\label{moser-besov}Let $s>0$ and $q\in[1,\infty]$.  Then we have
\begin{equation*}\|fg\|_{B^s_{p,q}}\leq C(\|f\|_{L^{p_1}}\|g\|_{B^s_{p_2,q}}
+\|g\|_{L^{r_1}}\|f\|_{B^s_{r_2,q}})
\end{equation*}
where $p_i,r_i\in[1,\infty]$ and $1/p=1/p_1+1/p_2=1/r_1+1/r_2$.
\end{proposition}

Using Proposition $\ref{moser-besov}$ and ($\ref{besov
embedding}$), we have, for $1/p=1/p_1+1/p_2$ and $p_1\leq p_2$,
\begin{equation*}\aligned\|u^2\|_{B^s_{p,q}}&\leq C\|u\|_{L^{p_2}}\|u\|_{B^s_{p_1,q}}
\\&\leq C\|u\|_{B^{\varepsilon}_{p_2,q}}\|u\|_{B^s_{p_1,q}}
\\&\leq{C}\|u\|^2_{B^s_{p_1,q}}
\endaligned
\end{equation*}
where $\varepsilon>0$ and $s=\varepsilon+n(2/p_1-1/p)$.  We record this as a proposition.
\begin{proposition}\label{besovproduct}Let $s>0$ and $p,p_1,q\in[1,\infty]$.  Then
\begin{equation*}\|u^2\|_{B^s_{p,q}}\leq C\|u\|^2_{B^s_{p_1,q}}
\end{equation*}
where $p<p_1\leq 2p$ and $s>n(2/p_1-1/p)$.
\end{proposition}

Next we state two Bernstein inequalities (see Appendix $A$ in \cite{taonde}).  We let $\Lambda=\sqrt{-\triangle}$, $\alpha\geq 0$, and $1\leq p\leq q\leq\infty$.  If $\text{supp}~\hat{f}\subset\{\xi\in\mathbb{R}^n:|\xi|\leq
K2^j\}$ and $\text{supp}~\hat{g}\subset\{\xi\in\mathbb{R}^n:K_12^j\leq |\xi|\leq
K_22^j\}$ for some $K, K_1, K_2>0$ and some integer $j$, then
\begin{equation}\label{Bernstein}\aligned \|\Lambda^\alpha f\|_q\leq
C2^{j\alpha+jn(1/p-1/q)}\|f\|_p
\\ \tilde{C}2^{j\alpha+jn(1/p-1/q)}\|f\|_p\leq \|\Lambda^\alpha f\|_q\leq
C2^{j\alpha+jn(1/p-1/q)}\|f\|_p.
\endaligned 
\end{equation}

We conclude with some results governing the behavior of the heat kernel on Besov spaces.

\begin{proposition}\label{hkb} Let $1\leq p_1\leq p_2<\infty$,
$-\infty<s_1\leq s_2<\infty$, and let $0<q<\infty$.  Then
\begin{equation*}\|e^{t\triangle} u\|_{B^{s_2}_{p_2,q}}\leq
Ct^{-(s_2-s_1+n/p_1-n/p_2)/2}\|u\|_{B^{s_1}_{p_1,q}}.
\end{equation*}
\end{proposition}

Using the Sobolev space heat kernel estimate, we get 
\begin{equation*}\aligned \|e^{t\triangle}u\|_{B^{s_2}_{p_2,q}}
&=\|e^{t\triangle} u\|_{L^{p_2}}+\(\sum
(2^{js_1}\|2^{j(s_2-s_1)}\phi_j *
e^{t\triangle}u\|_{L^{p_2}})^q\)^{1/q}
\\ &\leq t^{(n/p_1-n/p_2)/2}\|u\|_{L^{p_1}}+\(\sum
(2^{js_1}\|\phi_j * e^{t\triangle} u\|_{H^{s_2-s_1,p_2}})^q\)^{1/q}
\\ &\leq t^{-(n/p_1-n/p_2)/2}\|\phi * u\|_{L^{p_1}}+t^{\sigma}\(\sum
(2^{js_1}\|\phi_j * u\|_{L^{p_1}})^q\)^{1/q}
\\ &\leq
t^{\sigma}\|u\|_{B^{s_1}_{p_1,q}}.
\endaligned
\end{equation*}
where $\sigma=-(s_2-s_1+n/p_1-n/p_2)/2$.

We have an immediate corollary.
\begin{corollary}\label{hkb cor}With the parameters as in Proposition $\ref{hkb}$, we have
\begin{equation*}\lim_{t\rightarrow 0}t^{\gamma/2}\|e^{t\triangle}f\|_{\dot{B}^{\alpha+\beta}_{r,q_1}}=0
\end{equation*}
provided $\sigma>0$.
\end{corollary}

\section{Proof of Theorem $\ref{main global piece}$}\label{global extension}
In this section we establish an $\emph{a~priori}$ Besov space bound on solutions to the LANS equation.  This argument requires qualitatively different methods than those used in \cite{MS} and \cite{sobpaper} to obtain comparable bounds in the Sobolev space norm.  Moreover, the result applies in dimension $n\geq 3$, an improvement over the requirement that $n=3$ for Sobolev spaces.  

We begin by recalling Bony's notion of paraproduct (see $0.17$ in Chapter
2 of \cite{TT}). We have that $fg=T_fg+T_gf+R(f,g)$, where
\begin{equation}\label{paraproduct}T_f g=\sum_k (S_{k-2}f)\triangle_{k}g
\end{equation}
and
\begin{equation*}R(f,g)=\sum_k (\sum_{l=k-1}^{k+1}\triangle_l f)(\triangle_k g).
\end{equation*}

The following calculation will be of use throughout this section.  Using Bony's paraproduct and equations $(\ref{besovlemma1})$ and $(\ref{besovpieces67})$, we have that
\begin{equation*}\aligned\triangle_j(fg)&= \sum_{|j-k|\leq
2}\triangle_j(S_{k-2}f\triangle_{k} g)+\sum_{|j-k|\leq
2}\triangle_j(S_{k-2}g\triangle_{k} f)
\\&+\sum_{k\geq j-3}\triangle_j(\triangle_k
g\sum_{l=k-1}^{k+1}\triangle_{l} f)
\\&=I+II+III.
\endaligned
\end{equation*}
Applying Young's inequality and then Holder's inequality, we have
\begin{equation}\label{bony1}\aligned\|I\|_p&\leq C\sum_{|j-k|\leq 2}\|\psi_j\|_1\|S_{k-2}f\triangle_k
g\|_p
\\ &\leq C\sum_{|j-k|\leq 2}\|S_{k-2}f\|_\infty\|\triangle_k
g\|_p,
\endaligned
\end{equation}
and similarly
\begin{equation}\label{bony2}\|II\|_p\leq C\sum_{|j-k|\leq
2}\|S_{k-2}g\|_\infty\|\triangle_k f\|_p.
\end{equation}
\begin{equation}\label{bony3}\|III\|_p\leq C\sum_{k\geq j-3}
(\sum_{l=k-1}^{k+1}\|\triangle_k g\|_p\|\triangle_l f\|_\infty).
\end{equation}

We pause here to comment on our strategy.  To prove Theorem $\ref{main global piece}$, we need bounds on the two separate pieces of the Besov space norm.  We will begin by working with the more complicated $\|\cdot\|_{\tilde{B}^s_{p,q}}$ piece.   

The following proposition is the central piece to our desired estimate.  
\begin{proposition}\label{besovpropglobal}Any solution $u$ to $(\ref{LANS})$ satisfies
\begin{equation*}\frac{d}{dt}\|u\|^q_{\tilde{B}^r_{2,q}}\leq
C\|u\|_{\tilde{B}^{1+n/2}_{2,q}}\|u\|^q_{\tilde{B}^r_{2,q}},
\end{equation*}
provided $r>2$.
\end{proposition}

The proof of the Proposition is the content of the next Section.

\subsection{Proof of Proposition $\ref{besovpropglobal}$}
We begin with an equivalent form of the
LANS equations (see Section $3$ of \cite{MS}):
\begin{equation*}\aligned &\partial_t(1-\alpha^2\triangle)u
+\nabla_u (1-\alpha^2\triangle)u-\alpha^2(\nabla u)^T\cdot\triangle
u
\\ =&-\nu(1-\alpha^2\triangle)P^\alpha \triangle u-\nabla p,
\endaligned
\end{equation*}

Before proceeding, we recall the definitions of some of these objects.  First, $\nabla_u v=(v\cdot \nabla)u$ is the vector with $j^{th}$ component $\sum_{i=1}^n v_i\partial_{x_i} u_j$.  Secondly, $(\nabla u)^T \cdot v$ is the vector with $j^{th}$ component equal to the $\mathbb{R}^n$ dot product of the vectors $\partial_j u$ and $v$.

Now, applying $\triangle_j$ to both sides and taking the $L^2$ inner
product with $2\triangle_j u$, we get
\begin{equation}\label{global1}I_1+I_2+I_3+I_4=I_5
\end{equation}
where
\begin{equation*}\aligned I_1&=(\partial_t(1-\alpha^2\triangle)\triangle_j u,\triangle_j u)_{L^2},
\\ I_2&=(\triangle_j(\nabla_u(1-\alpha^2\triangle)u,\triangle_j u)_{L^2},
\\ I_3&=-\alpha^2(\triangle_j((\nabla u)^T\cdot\triangle u,\triangle_j u)_{L^2},
\\ I_4&=\nu((1-\alpha^2\triangle)P^\alpha\triangle \triangle_j u,\triangle_j u)_{L^2},
\\ I_5&=-(\nabla \triangle_j p,\triangle_j u)_{L^2}.
\endaligned
\end{equation*}

Applying integration by parts to $I_1$, $I_4$, and $I_5$, we get 
\begin{equation}\label{global2}\aligned&\frac{1}{2}\frac{d}{dt}(\|\triangle_j u\|^2_{L^2}
+\alpha^2\|A^{1/2}\triangle_j u\|^2_{L^2}) \leq &|I_2|+|I_3|.
\endaligned
\end{equation}
To proceed, we need to estimate $I_2$ and $I_3$.

\subsection*{Step 1: $I_2$ and $I_3$ estimates}\label{$I_2$ and $I_3$
estimates} We begin with $I_2$, and we have by Holder's inequality
that
\begin{equation*}|I_2|\leq C\|\triangle_j u\|_2
\|\triangle_j (\nabla_u(1-\alpha^2\triangle)u)\|_2.
\end{equation*}
To estimate the second term, we use $(\ref{bony1})$,
$(\ref{bony2})$, and $(\ref{bony3})$ to write
\begin{equation*}\triangle_j (\nabla_u(1-\alpha^2\triangle)u)=J_1+J_2+J_3
\end{equation*}
where
\begin{equation*}\aligned J_1&= \sum_{|j-k|\leq 2}\triangle_j((S_{k-2}u\cdot \nabla)
\triangle_k u) -\alpha^2\sum_{|j-k|\leq
2}\triangle_j((S_{k-2}u\cdot \nabla)\triangle_k \triangle u)
\\J_2&=\sum_{|j-k\leq 2}\triangle_j ((\triangle_k u\cdot\nabla)S_{k-2}u)-\alpha^2
\sum_{|j-k|\leq 2}\triangle_j ((\triangle_k u\cdot\nabla)S_{k-2}\triangle u)
\\J_3&=\sum_{k\geq j-3}\triangle_j \(\triangle_k u
\sum_{l=-1}^{l=1}\nabla\triangle_{k+l} u\) -\alpha^2\sum_{k\geq
j-3}\triangle_j \(\triangle_k u
\sum_{l=-1}^{l=1}\nabla\triangle_{k+l} \triangle u\).
\endaligned
\end{equation*}

For notational convenience, we define $J_{i,j}$, with $j=1,2$, to be
the $j^{th}$ of the two summations in $J_i$.  We begin with
$J_{1,2}$. Recalling that $\psi_j$ is the convolution kernel for
$\triangle_j$, we have by integration by parts and the
incompressibility condition that
\begin{equation*}\aligned &\int\psi_j(x-y)((S_{k-2}u(y))\nabla\(\triangle\triangle_k
 u(y)\)dy
\\ =&\int
-\triangle(\nabla\psi_j(x-y)\cdot(S_{k-2}u(y)))\triangle_k
u(y)dy.
\endaligned
\end{equation*}
Next, we use the product rule to distribute the Laplacian through
the product, apply Young's inequality, and then take the $L^\infty$
norm of the pieces involving $S_{k-2} u$ and its derivatives.
Recalling that the $L^1$ norm of $\psi_j$ and its derivatives are
independent of $j$, we have
\begin{equation*}\aligned &\|\sum_{|j-k|\leq 2}\triangle_j((S_{k-2}u\cdot
\nabla)\triangle_k \triangle u)\|_{L^2}
\\ \leq &C\sum_{|j-k|\leq 2}(\|S_{k-2} u\|_{L^\infty}+\|\nabla S_{k-2} u\|_{L^\infty}
+\|\triangle S_{k-2} u\|_{L^{\infty}})\|\triangle_k u\|_{L^2}.
\endaligned
\end{equation*}

Recalling the definition of $S_{k-2}$ and using Bernstein's
inequality (equation (\ref{Bernstein})) we get
\begin{equation*}\aligned &\|\sum_{|j-k|\leq 2}\triangle_j((S_{k-2}u\cdot
\nabla)\triangle_k \triangle u)\|_{L^2}
\\ \leq &C\sum_{|j-k|\leq 2}\|\triangle_k
u\|_{L^2}\sum_{m<k-2}2^{(2+n/2)m}\|\triangle_m u\|_{L^2}.
\endaligned
\end{equation*}

For $J_{1,1}$, a similar computation gives
\begin{equation*}\aligned &\|\triangle_j((S_{k-2}u\cdot \nabla)\triangle_k u)\|_{L^2}
\\ \leq &C\|\triangle_k u\|_{L^2}\sum_{m<k-2}2^{nm/2}\|\triangle_m
u\|_{L^2}.
\endaligned
\end{equation*}
So we finally get that $J_1$ satisfies
\begin{equation*}|J_1|\leq C\sum_{|j-k|\leq 2}\|\triangle_k u\|_{L^2}
\sum_{m<k-2}2^{(2+n/2)m}\|\triangle_m u\|_{L^2}.
\end{equation*}

$J_2$ satisfies the same estimate, so we turn to $J_3$.  Again using
integration by parts, Young's inequality, and then Bernstein's inequality, we have
\begin{equation*}\|J_{3,2}\|_{L^2}\leq C\sum_{k\geq j-3}
\|\triangle_k
u\|_{L^2}\sum_{l=-1}^{l=1}2^{(2+n/p)(k+l)}\|\triangle_{k+l}
u\|_{L^2}
\end{equation*}
and 
\begin{equation*}\|J_{3,1}\|_{L^2}
\leq C\sum_{k\geq j-3}\|\triangle_k
u\|_{L^2}\sum_{l=-1}^{l=1}2^{(n/p)(k+l)j}\|\triangle_{k+l}
u\|_{L^2}.
\end{equation*}

So $J_3$ satisfies
\begin{equation*}|J_{3}|\leq C\sum_{k\geq j-3}
\|\triangle_k
u\|_{L^2}\sum_{l=-1}^{l=1}2^{(2+n/p)(k+l)}\|\triangle_{k+l}
u\|_{L^2}.
\end{equation*}

So we finally estimate $I_2$ by
\begin{equation} \label{I2}\aligned |I_2|&=C\|\triangle_j u\|_{L^2}\|J_1+J_2+J_3\|_{L^2}
\\ &\leq C\|\triangle_j u\|_{L^2}\sum_{|j-k|\leq 2}\|\triangle_k u\|_{L^2}
\sum_{m<k-2}2^{(2+n/2)m}\|\triangle_m u\|_{L^2}
\\+&C\|\triangle_j
u\|_{L^2}\sum_{k\geq j-3} \|\triangle_k
u\|_{L^2}\sum_{l=-1}^{l=1}2^{(2+n/p)(k+l)} \|\triangle_{k+l}
u\|_{L^2}.
\endaligned
\end{equation}

The estimation of $I_3$ is similar, so the details will be omitted.
The key difference between the two estimates is that, in the case of
$I_2$, integrating the gradient term by parts and applying the
incompressibility condition essentially removed one of the three
derivatives from $I_2$. Because the gradient term in $I_3$ is
$(\nabla u)^T$ instead of $\nabla u$, this does not work when
estimating $I_3$, and the result is the presence of an extra
derivative.  The estimate for $I_3$ is
\begin{equation} \label{I3}\aligned |I_3|&\leq C\|\triangle_j u\|_{L^2}
\sum_{|j-k|\leq 2}\|\triangle_k u\|_{L^2}
\sum_{m<k-2}2^{(3+n/2)m}\|\triangle_m u\|_{L^2}
\\+&C\|\triangle_j u\|_{L^2}\sum_{k\geq j-3} \|\triangle_k
u\|_{L^2}\sum_{l=-1}^{l=1}2^{(3+n/p)(k+l)} \|\triangle_{k+l}
u\|_{L^2}.
\endaligned
\end{equation}

Using $(\ref{I2})$ and $(\ref{I3})$ in $(\ref{global2})$ we get
\begin{equation}\label{bernlike}\aligned &\frac{d}{dt}\(\|\triangle_j u\|^2_{L^2}
+\alpha^2\|A^{1/2}\triangle_j u\|_{L^2}^2\)
\\ &\leq C\|\triangle_j u\|_{L^2}
\sum_{|j-k|\leq 2}\|\triangle_k u\|_{L^2}
\sum_{m<k-2}2^{(3+n/2)m}\|\triangle_m u\|_{L^2}
\\+&C\|\triangle_j u\|_{L^2}\sum_{k\geq j-3} \|\triangle_k
u\|_{L^2}\sum_{l=-1}^{l=1}2^{(3+n/p)(k+l)} \|\triangle_{k+l}
u\|_{L^2}.
\endaligned
\end{equation}
In the next section, we work on the left-hand side of
$(\ref{bernlike})$.

\subsection*{Step 2: Exploiting the LANS term}\label{Exploiting the LANS term}
By Plancherel's Theorem $(\ref{bernlike})$ becomes 
\begin{equation*}\aligned 2^{2j}\frac{d}{dt}\(\|\triangle_j u\|^2_{L^2}\)
&\leq C\|\triangle_j u\|_{L^2} \sum_{|j-k|\leq 2}\|\triangle_k
u\|_{L^2} \sum_{m<k-2}2^{(3+n/2)m}\|\triangle_m u\|_{L^2}
\\+&C\|\triangle_j u\|_{L^2}\sum_{k\geq j-3} \|\triangle_k
u\|_{L^2}\sum_{l=-1}^{l=1}2^{(3+n/2)(k+l)} \|\triangle_{k+l}
u\|_{L^2}.
\endaligned
\end{equation*}

Computing the time derivative, multiplying both sides by $q2^{rjq}\|\triangle_j
u\|^{q-1}_{L^2}$ and summing over $j\geq 0$,  we get
\begin{equation}\label{globalest}\frac{d}{dt}\|u\|^q_{\tilde{B}^r_{2,q}}\leq
K_1+K_2,
\end{equation}
where
\begin{equation*}K_1=Cq\sum_{j\geq 0} 2^{rjq}2^{-2j}\|\triangle_j u\|^{q-1}_{L^2}
\sum_{|j-k|\leq 2}\|\triangle_k u\|_{L^2}
\sum_{m<k-2}2^{(3+n/2)m}\|\triangle_m u\|_{L^2},
\end{equation*}
and
\begin{equation*}K_2=Cq\sum_{j\geq 0} 2^{rjq}2^{-2j}\|\triangle_j
u\|^{q-1}_{L^2}\sum_{k\geq j-3} \|\triangle_k
u\|_{L^2}\sum_{l=-1}^{l=1}2^{(3+n/p)(k+l)} \|\triangle_{k+l}
u\|_{L^2}.
\end{equation*}

We remark that the $2^{-2j}$ term in $K_1$ and $K_2$ is the ``gain" we get by using the LANS equations instead of the
Navier-Stokes equations.

\subsection*{Step 3: Estimating $K_1$ and $K_2$}\label{Estimating $K_1$ and
$K_2$}

We begin by re-writing $K_1$ as
\begin{equation*}\aligned K_1=Cq\sum_{j\geq 0} 2^{rjq}\|\triangle_j u\|^{q-1}_{L^2}
\sum_{k=-2}^2\|\triangle_{j+k} u\|_{L^2}
\sum_{m<j+k-2}2^{(3+n/2)m-2j}\|\triangle_m u\|_{L^2}.
\endaligned
\end{equation*}

We start by working on the last summation.  We have by Holder's
inequality that
\begin{equation*}\aligned &\sum_{m<j+k-2}2^{(1+n/2)m}2^{2m-2j}\|\triangle_m
u\|_2
\\ &\leq C\(\sum_{m<j+k-2}2^{(1+n/2)qm}\|\triangle_m
u\|^q_2\)^{1/q}\(\sum_{m<j+k-2}2^{2q'(m-(j+k-2)}2^{2(k-2)q'}\)^{1/q'}
\\ &\leq C\|u\|_{\tilde{B}^{1+n/2}_{2,q}},
\endaligned
\end{equation*}
where $q'$ is the Holder conjugate exponent to $q$.

Returning to $K_1$, we have
\begin{equation*}|K_1|\leq C\|u\|_{\tilde{B}^{1+n/2}_{2,q}}
\sum_{j\geq 0} 2^{rjq}\|\triangle_j u\|^q_{L^2}
\sum_{k=-2}^2\|\triangle_{j+k} u\|_{L^2}.
\end{equation*}

For the remaining summation, we have
\begin{equation*}\aligned &\sum_{j\geq 0} 2^{rjq}\|\triangle_j u\|^{q-1}_{L^2}
\sum_{k=-2}^2\|\triangle_{j+k} u\|_{L^2}
\\ &\leq C\(\sum_{j\geq 0} 2^{rj(q-1)q'}\|\triangle_j
u\|^{q'(q-1)}_{L^2}\)^{q'}\(\sum_j
(\sum_{k=-2}^2 2^{r(j+k)}2^{-rk}\|\triangle_{j+k}
u\|_{L^2})^q\)^q
\\ &\leq C\|u\|_{\tilde{B}^r_{2,q}}^{q-1}\|u\|_{\tilde{B^r_{2,q}}}
\leq C\|u\|^q_{\tilde{B}^r_{2,q}},
\endaligned
\end{equation*}
where $q'$ is again the Holder conjugate exponent to $q$.  So we
finally bound $K_1$ by
\begin{equation*}K_1\leq
C\|u\|^q_{\tilde{B}^r_{2,q}}\|u\|_{\tilde{B}^{1+n/2}_{2,q}}.
\end{equation*}

Bounding $K_2$ the same way, we have
\begin{equation}\label{reason for trouble}
|K_2|\leq C\|u\|_{B^r_{2,q}}^{q-1} \tilde{K},
\end{equation}
where
\begin{equation*}\tilde{K}=\(\sum_j \(\sum_{k>-3}\sum_{l=-1}^1 2^{(r-2)j}2^{(3+n/2)(j+k+l)}
\|\triangle_{j+k}u \|_{L^2}\|\triangle_{j+k+l}
u\|_{L^2}\)^q\)^{1/q}.
\end{equation*}

A straightforward calculation shows that $\tilde{K}$ satisfies 
\begin{equation*}\aligned \tilde{K}\leq C\|u\|_{B^{1+n/2}_{2,q}} \(\sum_j \(\sum_{k>-3}2^{k(2-r)}\sum_{l=-1}^1
(2^{r(j+k)}\|\triangle_{j+k} u\|_2)\)^q\)^{1/q}.
\endaligned
\end{equation*}

Applying Minkowski's inequality, we get
\begin{equation}\label{cmpunkorton}\aligned |\tilde{K}|\leq C\|u\|_{\tilde{B}^{1+n/2}_{2,q}}\|u\|_{\tilde{B}^r_{2,q}}\sum_{k>-3} 2^{k(2-r)}.
\endaligned
\end{equation}
This last sum will be finite provided $r>2$.  We remark that this 
restriction on $r$ could be lifted by allowing for additional
regularity in the $B^{1+n/2}_{2,q}$ term.  Specifically, by choosing
$\varepsilon$ such that $2<r+\varepsilon$, we would have
\begin{equation*}|\tilde{K}|\leq
\|u\|_{\tilde{B}^{1+n/2+\varepsilon}_{2,q}}\|u\|_{\tilde{B}^r_{2,q}}.
\end{equation*}
Using this estimate would provide a slightly different result for Proposition $\ref{besovpropglobal}$.

Assuming $r>2$, we plug these estimates back into the $K_2$ estimate, and get
\begin{equation*}|K_2|\leq
C\|u\|^q_{\tilde{B}^r_{2,q}}\|u\|_{\tilde{B}^{1+n/p}_{2,q}}
\end{equation*}
provided $r>2$.

Plugging the $K_1$ and $K_2$ estimates into $(\ref{globalest})$, we
finally get
\begin{equation}\label{keyglobalpiece}\frac{d}{dt}\|u\|^q_{\tilde{B}^r_{2,q}}\leq
C\|u\|^q_{\tilde{B}^{r}_{2,q}}\|u\|_{\tilde{B}^{1+n/2}_{2,q}}
\end{equation}
for $r>2$.

\subsection{Proof of Theorem $\ref{main global piece}$}\label{Proof of Global
Existence} Now we return to proving Theorem $\ref{main global piece}$.  Recalling our construction of Besov spaces, we need to show that 
\begin{equation}\label{sammy}\|\Psi\ast u(t)\|_{L^2}\leq \|u_0\|_{B^r_{2,q}}\exp \(C\int_I \|u(t)\|_{B^{1+n/2}_{2,q}} dt\)
\end{equation}
and 
\begin{equation}\label{bil}\|u(t)\|_{\tilde{B}^{n/2}_{2,q}}\leq C\|u_0\|_{B^r_{2,q}}\exp \(C\int_I \|u(t)\|_{B^{1+n/2}_{2,q}} dt\).
\end{equation}

To prove $(\ref{bil})$, we re-write $(\ref{keyglobalpiece})$ as
\begin{equation}\label{keyglobalpiece2}  \frac{d}{dt}\|u\|_{\tilde{B}^r_{2,q}}\leq
C\|u\|_{\tilde{B}^{r}_{2,q}}\|u\|_{\tilde{B}^{1+n/2}_{2,q}}
\end{equation}
  
Applying Gronwall's inequality to $(\ref{keyglobalpiece2})$, we get
\begin{equation}\aligned\label{globalbesovhompiece} \|u(t)\|_{\tilde{B}^{r}_{2,q}}
&\leq \|u_0\|_{\tilde{B}^{r}_{2,q}}
\exp\(C\int_0^T\|u(s)\|_{\tilde{B}^{1+n/2}_{2,q}}ds\)
\\ &\leq \|u_0\|_{{B}^{r}_{2,q}}
\exp\(C\int_0^T\|u(s)\|_{{B}^{1+n/2}_{2,q}}ds\)
\endaligned
\end{equation}
which proves $(\ref{bil})$.

To prove $(\ref{sammy})$, we go back to the LANS equation (see Section $3$ of \cite{MS}):
\begin{equation}\label{LANSglobalpi}\aligned &\partial_t (1-\alpha^2\triangle)u
+\nabla_u[(1-\alpha^2\triangle)u] -\alpha^2(\nabla
u)^T\cdot\triangle u
\\=&-(1-\alpha^2\triangle)A u-\nabla p.
\endaligned
\end{equation}

%start of H^1 bound
Taking the $L^2$ product of $(\ref{LANSglobalpi})$ with
$u$, we get
\begin{equation*}I_1+I_2+I_3=J_1+J_2
\end{equation*}
where
\begin{equation*}\aligned I_1&=(\partial_t (1-\alpha^2\triangle)u,u)
\\ I_2&=(\nabla_u u,u)
\\ I_3&=-\alpha^2\((\nabla_u \triangle u,u)+((\nabla u)^T\cdot \triangle u, u)\)
\\ J_1&=-((1-\alpha^2 \triangle)(A u),u)
\\ J_2&=(\nabla p,u).
\endaligned
\end{equation*}

We start with $I_1$, which becomes
\begin{equation*}\aligned I_1&=(\partial_t u,u)-\alpha^2(\triangle \partial_t u,u)
\\ &=\frac{1}{2}\partial_t (\|u\|_{L^2}^2+\alpha^2\|A^{1/2}u\|^2_{L^2}).
\endaligned
\end{equation*}
Applying integration by parts to $I_2$, $I_3$, and $J_2$ and recalling that $\div u=0$, we get that all three terms vanish.  For $J_1$ we have 
\begin{equation*}J_1=-((1-\alpha^2 \triangle)(A u),u)=-(A^{1/2}u,A^{1/2}u)-\alpha^2(Au,Au).
\end{equation*}

Applying this to $(\ref{LANSglobalpi})$, we get
\begin{equation}\label{global11}\frac{1}{2}\partial_t (\|u(t)\|^2_{L^2}+\alpha^2\|u(t)\|^2_{\dot{H}^{1,2}})
\leq -(\|A^{1/2}u(t)\|^2_{L^2}+\alpha^2\|Au(t)\|_{L^2}^2),
\end{equation}
where $\dot{H}$ denotes the homogeneous Sobolev norm.  This proves
that $\|u(t)\|_{H^{1,2}}$ is decreasing in time.

Using $(\ref{global11})$, we have 
\begin{equation}\label{cmpunksammy}\|\Psi\ast u(t)\|_{L^2}\leq \|u(t)\|_{H^{1,2}}\leq M,
\end{equation}
which is stronger than $(\ref{sammy})$.  This proves Theorem $\ref{main global piece}$.

\section{Local solutions in $\dot{C}^T_{a;s,p,q}$}\label{local continuous solutions}
We begin by re-writing the LANS equation as 
\begin{equation}\label{LANS2}\aligned\partial_t u-Au+P^\alpha (\div\cdot(u\tensor
u)+\div\tau^\alpha u)=0,
\endaligned
\end{equation}
where the recurring terms are as in $(\ref{LANS})$, with the exception that we set $\nu=1$.  For the new terms, we set $A=P^\alpha\triangle$, $u\tensor u$ is the tensor with
$jk$-components $u_ju_k$ and $\text{div}\cdot(u\tensor u)$ is the
vector with $j$-component $\sum_k\partial_k(u_ju_k)$.  $P^\alpha$ is
the Stokes Projector, defined as
\begin{equation*}P^\alpha(w)=w-(1-\alpha^2\triangle)^{-1}\nabla f
\end{equation*}
where $f$ is a solution of the Stokes problem: Given $w$, there is a
unique divergence-free $v$ and a unique (up to additive constants) function $f$ such
that
\begin{equation*}(1-\alpha^2\triangle)v+\nabla f=(1-\alpha^2\triangle)w.
\end{equation*}
For a more explicit treatment of the Stokes
Projector, see Theorem 4 of \cite{SK}.

Using Duhamel's principle, we write ($\ref{LANS2}$) as the
integral equation
\begin{equation}\label{intversion}u=\Gamma\varphi-G\cdot P^\alpha(\div(u\tensor u+\tau^\alpha(u))) 
\end{equation}
with
\begin{equation*}(\Gamma\varphi)(t)=e^{tA}\varphi,
\end{equation*}
where $A$ agrees with $\triangle$ when restricted to $P^\alpha
H^{r,p}$, and  
\begin{equation*}G\cdot g(t)=\int_0^t e^{(t-s)A}\cdot
g(s)ds.
\end{equation*}

To prove Theorem $\ref{special case 1}$ and Theorem $\ref{thm1}$, we will apply the standard contraction method.  Section $\ref{Preliminary Results}$ contains the supporting details necessary for the proof, which can be found in Section $\ref{Proof of Theorem 1besov}$.  These results follow the same outline as that from \cite{sobpaper}, and more detailed arguments can be found there.     

\subsection{Preliminary Estimates}\label{Preliminary Results}

Our first calculation is a result for $\tau^\alpha$.
\begin{lemma}\label{tau lemma besov}Given $r>1$, $1\leq q<\infty$ and
$p,\bar{p}\in (1,\infty)$, where
\begin{equation*}\aligned p&\leq 2\bar{p},
\\ \bar{s}&:=n(2/p-1/\bar{p}),
\\ 0&\leq \bar{s}< r-1,
\endaligned
\end{equation*}
we have $\div\tau^\alpha:B^r_{p,q}\rightarrow B^r_{\bar{p},q}$.
Specifically, we have the estimate
\begin{equation*}\|\div\tau^\alpha(u)\|_{B^r_{\bar{p},q}}\leq C\|u\|^2_{B^r_{p,q}}.
\end{equation*}
\end{lemma}

We have by Proposition $\ref{besovproduct}$ that
\begin{equation*}\aligned\|\div \tau^\alpha(u)\|_{B^r_{\bar{p},q}}&\leq C\|\tau^\alpha(u)\|_{B^{r+1}_{\bar{p},q}}
\\ &\leq C\|Def(u)\cdot Rot(u)\|_{B^{r-1}_{\bar{p},q}}
\\ &\leq C\|\nabla u\|^2_{B^{r-1}_{p,q}} \leq C\|u\|^2_{B^{r}_{p,q}}.
\endaligned
\end{equation*}

We remark that this differs from the Sobolev result in that $r>1$ instead of $r\geq 1$ and that $\bar{s}< r-1$ instead of $\bar{s}\leq r-1$.  This Lemma has an immediate corollary.
\begin{corollary}\label{tau cor besov}$\div \tau^\alpha :\dot{C}^T_{a;r,p,q}\rightarrow \dot{C}^T_{2a;r,\bar{p},q}$,
with the estimate
\begin{equation*}\|\div \tau^\alpha(u)\|_{2a;r,\bar{p},q}\leq
C\|u\|^2_{a;r,p,q}.
\end{equation*}
\end{corollary}

Next, we record some results for the operator $V^\alpha$.
\begin{proposition}\label{Vprop besov}With the parameters as in Lemma $\ref{tau lemma besov}$, we have 
\begin{equation*}V^\alpha:\dot{C}^T_{a;r,p,q}\times\dot{C}^T_{a;r,p,q}
\rightarrow \dot{C}^T_{2a;r-1,p',q}
\end{equation*}
with the estimate
\begin{equation}\label{Vest besov}\|V^\alpha(u,v)\|_{2a;s-1,p',q}
\leq \|u\|_{a;s,p,q}\|v\|_{a;s,p,q}.
\end{equation}
\end{proposition}
This follows from a calculation parallel to the one used to prove
Proposition $1$ in \cite{sobpaper} with $(\ref{besov embedding})$
replacing the Sobolev embedding results.

\begin{corollary}\label{Vcor besov}With the same assumptions
on the parameters as in Proposition $\ref{Vprop besov}$, we have
that
\begin{equation}\label{Vcont besov}\aligned
&\|V^\alpha(u(s))-V^\alpha(v(s))\|_{a;s-1,p',q}
\\ \leq{C}&(\|u\|_{a/2;s,p,q}+\|v\|_{a/2;s,p,q})\|u-v\|_{a/2;s,p,q}.
\endaligned
\end{equation}
\end{corollary}

For the remainder of the section, we impose the following list of
restrictions on our parameters.  We have
\begin{equation}\label{hkb conditions}\aligned-\infty&<s_0\leq s_1<\infty
\\ 1&\leq q\leq \infty
\\ 1&\leq p_0\leq p_1<\infty
\\ \sigma&= s_1-s_0+n(1/p_0-1/p_1).
\endaligned
\end{equation}

The following proposition is an immediate consequence of Corollary
$\ref{hkb cor}$.
\begin{proposition}\label{gamma prop besov}Provided the parameters satisfy
$(\ref{hkb conditions})$, we have that $\Gamma
:B^{s_0}_{p_0,q}\rightarrow \dot{C}^T_{\sigma/2;s_1,p_1,q}$.
\end{proposition}

Lastly, we turn our attention to the operator $G$. Using Proposition $\ref{hkb}$, we have
\begin{equation*}\aligned\|G\cdot u\|_{B^{s_1}_{p_1,q}}
&\leq C\int_0^t (t-s)^{-\sigma/2}\|u\|_{B^{s_0}_{p_0,q}}ds
\\&\leq C\|u\|_{k;s_0,p_0,q}\int_0^t(t-s)^{-\sigma/2}t^{-k}
\\&\leq Ct^{-\sigma/2-k+1}\|u\|_{k;s_0,p_0,q}.
\endaligned
\end{equation*}
We record this as a proposition.
\begin{proposition}\label{G Prop Besov}With parameters as specified in $(\ref{hkb conditions})$,
$0<\sigma/2<1$, $0\leq k_0<1$, and $k_1=k_0+\sigma/2-1$, we have
\begin{equation*}\|G\cdot u\|_{k_1;s_1,p_1,q}\leq C\|u\|_{k_0;s_0,p_0,q}.
\end{equation*}
\end{proposition}

\subsection{Proof of Theorems $\ref{special case 1}$ and $\ref{thm1}$}\label{Proof of
Theorem 1besov}
As usual, we begin with the nonlinear map
\begin{equation*}\Phi u=\Gamma \varphi-G\cdot P^\alpha(V^\alpha(u)),
\end{equation*}
initial data $u_0\in B^r_{p,q}$ and define
\begin{equation*}E_{T,M}=\{v\in \bar{C}_{r,p,q} \cap \dot{C}^T_{a;s,\tilde{p},q}
:\|v-\Gamma\varphi\|_{0;r,p,q}+\|v\|_{a;s,\tilde{p},q}\leq M\}.
\end{equation*}

The proof is a standard contraction mapping argument.  For the details, see \cite{sobpaper}.  We get that $\Phi$ will be a contraction on $E_{T,M}$ provided there exists a $\bar{b}>0$ such that the following list of conditions holds (where, for notational convenience, we set $\bar{s}=s-2+\bar{b}-r+n/p$):

\begin{equation*}\label{}\aligned
&1< p\leq \tilde{p}
\\ &1\leq q\leq\infty
\\&s>1,~~ \bar{b}\geq 1, 
\\&\bar{s}\tilde{p}<n
\\&0<2a=s-n/\tilde{p}-r+n/p<1
\\&0\leq \bar{s}< s-1
\\&1<\frac{n\tilde{p}}{2n-\bar{s}\tilde{p}}<\infty
\\& 1> \bar{b}-r+n/p
\\& 1\leq \bar{b}+\frac{n}{\tilde{p}}-\bar{s}<2
\\& 2-2\bar{b}+\bar{s}\leq\frac{n}{p}\leq 2-\bar{b}+\bar{s}.
\endaligned
\end{equation*}
We remark that the requirements $s>1$, $\bar{s}<s-1$, and $1>\bar{b}-r+n/p$ differ from the analogous Sobolev space requirements, which were $s\geq 1$, $\bar{s}\leq s-1$, and $1\geq \bar{b}-r+n/p$ respectively.  These relations combine to force $r>n/p$.  

This result is actually stronger than required for Theorem $\ref{thm1}$, but requires the presence of the artificial parameter $\bar{b}$ (though this is somewhat compensated for by the explicit formula for $a$).  To get the result of Theorem $\ref{thm1}$, we fix $\bar{b}=1$, which provides the minimal allowable $r$ (and we again set $\bar{s}=s-2+\bar{b}-r+n/p$): 
\begin{equation}\label{conditionsbesov1}\aligned
&1< p\leq \tilde{p}
\\&1\leq q\leq\infty
\\&s>1,~~ \bar{s}\tilde{p}<n
\\ &r>n/p
\\&0<2a=s-n/\tilde{p}-b<1
\\&0\leq \bar{s}< s-1
\\&1<\frac{n\tilde{p}}{2n-\bar{s}\tilde{p}c}<\infty
\\& 0\leq \frac{n}{\tilde{p}}-\bar{s}<1
\\& \bar{s}\leq\frac{n}{p}\leq 1+\bar{s}.
\endaligned
\end{equation}

As in Section $3$ of \cite{sobpaper}, this general method can be modified to obtain the results in Theorem $\ref{special case 1}$.  We also have the following result which extends the range of parameters allowed on the auxiliary space.

\begin{proposition}\label{higher regularity theorem}Let $u_0\in B^{s_1}_{p,q}(\mathbb{R}^n)$ be divergence-free, and let $s_1<1$.  Let $u$ be a solution to the LANS equation ($\ref{LANS})$ such that   
\begin{equation*}u\in BC([0,T):B^{s_1}_{p,q}(\mathbb{R}^n))\cap \dot{C}^T_{(s_2-s_1)/2;s_2,p,q},
\end{equation*}
where $0<s_2-s_1<1$ and $s_2\geq 1$.  Then for all $r\geq s_2$, we have that $u\in \dot{C}^T_{(r-s_1)/2;r,p,q}$.
\end{proposition}
We remark that the requirement $s_1<1$ can be removed and the assumption that $u\in BC([0,T):B^{s_1}_{p,q}(\mathbb{R}^n))\cap \dot{C}^T_{(s_2-s_1)/2;s_2,p,q}$ can be changed to $u\in BC([0,T):B^{s_1}_{\tilde{p},q}(\mathbb{R}^n))\cap \dot{C}^T_{a;s_2,\tilde{p},q}$ for $\tilde{p}>p$.

\section{Local solutions in $L^a((0,T):B^s_{p,q}(\mathbb{R}^n))$}\label{local integral solutions}

In this section we consider Theorems $\ref{special case 2}$ and $\ref{thm2}$.  Section $\ref{Preliminary results2}$ contains the necessary supporting estimates, and the theorems are proven in Section $\ref{proof of thm2}$.

\subsection{Preliminary Results}\label{Preliminary results2}In this section we establish integral-in-time results for Besov space.  The proofs are similar to those in \cite{sobpaper} used for the analogous operators.
\begin{proposition}\label{gamma prop 3}Let $1<p_0\leq p_1<\infty$,
$1\leq q<\infty$, $-\infty<s_0\leq s_1<\infty$, and assume
$0<(s_1-s_0+n/p_0-n/p_1)/2=1/\sigma$.  Then $\Gamma$ maps
$B^{s_0}_{p_0,q_0}$ continuously into
$L^\sigma((0,\infty):B^{s_1}_{p_1,q_1})$ with the estimate
\begin{equation*}\|\Gamma u\|_{L^\sigma((0,\infty):B^{s_1}_{p_1,q_1})}
\leq C\|u\|_{B^{s_0}_{p_0,q_0}}.
\end{equation*}
\end{proposition}

The proof is similar to Proposition $4$ in \cite{sobpaper}, with two main distinctions, both due to the differences in interpolation theory between Sobolev and Besov spaces.  The first is that we interpolate using $s_0$ instead of $p_0$.  The second difference is that we do not require $p_0\leq \sigma$, as we did in Proposition $4$ of \cite{sobpaper}.

Our next two results involve the operator $V^\alpha$.
\begin{proposition}\label{Vprop3}With the parameters $s, p, p'$ and $q$ as in Proposition $\ref{Vprop besov}$, we have
\begin{equation*}V^\alpha:L^\sigma((0,T):B^{s}_{p,q})\rightarrow L^{\sigma/2}((0,T):B^{s-1}_{p,q})
\end{equation*}
with the estimate
\begin{equation*}\(\int_0^T \|V^\alpha(u(s))\|^{\sigma/2}_{B^{s-1}_{p',q}}ds\)^{2/\sigma}\leq
\(\int_0^T\|u(s)\|^\sigma_{B^s_{p,q}} ds\)^{2/\sigma}.
\end{equation*}
\end{proposition}

\begin{corollary}If $u,v\in L^\sigma((0,T):B^s_{p,q})$, then
\begin{equation*}\aligned &\(\int_0^T\|V^\alpha(u(s))-V^\alpha(v(s))\|^{\sigma/2}_{B^{s-1}_{p,q}}ds\)^{2/\sigma}
\\ \leq &\(\int_0^T\|v(s)+u(s)\|^\sigma_{B^s_{p,q}} ds\)^{2/\sigma}\(\int_0^T\|v(s)-u(s)\|^\sigma _{B^s_{p,q}} ds\)^{2/\sigma}.
\endaligned
\end{equation*}
\end{corollary}

We conclude this section with estimates for $G$.
\begin{proposition}\label{Gprop3}Given $1\leq p_0\leq p_1<\infty$, $1\leq q<\infty$,
$-\infty<s_0\leq s_1<\infty$, $1<\sigma_0<\sigma_1<\infty$ and
$1/\sigma_0-1/\sigma_1=1-(s_1-s_0+n/p_0-n/p_1)/2$, for any $T\in
(0,\infty]$, $G$ sends $L^{\sigma_0}((0,T):B^{s_0}_{p_0,q_0})$ into
$L^{\sigma_1}((0,T):B^{s_1}_{p_1,q_1})$ with the estimate
\begin{equation*}\|G\cdot u\|_{L^{\sigma_1}((0,T):B^{s_1}_{p_1,q_1})}
\leq C\|u\|_{L^{\sigma_0}((0,T):B^{s_0}_{p_0,q_0})}.
\end{equation*}
\end{proposition}

\begin{proposition}\label{Gprop4}$1<p_0\leq p_1<\infty$, $1\leq q<\infty$,
$-\infty<s_0\leq s_1<\infty$, and assume $1/p_1\leq
1/\sigma=1-(s_1-s_0+n/p_0-n/p_1)/2=$.  Then $G$ maps
$L^\sigma((0,T):B^{s_0}_{p_0,q_0})$ continuously into
$BC([0,T):B^{s_1}_{p_1,q_1})$ with the estimate
\begin{equation*}\sup_{t\in[0,T)}\|G\cdot u(t)\|_{B^{s_1}_{p_1,q_1}}
\leq C\|u\|_{L^{\sigma}((0,T):B^{s_0}_{p_0,q_0})}.
\end{equation*}
\end{proposition}

\subsection{Proof of Theorems $\ref{special case 2}$ and $\ref{thm2}$}\label{proof of thm2}

As usual, we begin by defining a Banach space $X_{T,M}$ to be the
set of all $u\in BC([0,T):B^r_{p,q})\cap
L^\sigma((0,T):B^s_{\tilde{p},q})$ such that
\begin{equation*}\sup_t\|u(t)-\Gamma u_0\|_{B^r_{p,q}}
+\|u\|_{\sigma;s,\tilde{p},q}\leq M
\end{equation*}
and we define the operator $\Phi$ by
\begin{equation*}\Phi u(t)=\Gamma\varphi+G(V^\alpha (u(t)))
\end{equation*}
where $\varphi\in B^r_{p,q}$.

The required restrictions on the parameters are similar to those from Section $4$ of \cite{sobpaper}, and again we get that $\Phi$ will be a contraction provided there exists a $\bar{b}>0$ such that the following holds (for notational convenience, we set $\bar{s}=s-2+\bar{b}-r+n/p$):
We record these parameters:
\begin{equation*}\aligned& 1< p\leq \tilde{p}<\infty
\\&1\leq q\leq \infty
\\&s> 1,~~ \bar{b}\geq 1,~~\bar{s}\tilde{p}<n
\\&0<2/\sigma=s-n/\tilde{p}-b<1
\\&0\leq \bar{s}< s-1
\\&1<\frac{n\tilde{p}}{2n-\bar{s}\tilde{p}}<\infty
\\& 1> \bar{b}-r+n/p
\\& s-\bar{b}\leq \frac{n}{\tilde{p}}+r-n/p\leq s
\\& \sigma/2\leq p\leq \sigma.
\endaligned
\end{equation*}
As in Section $\ref{local continuous solutions}$, we get Theorem $\ref{thm2}$ by choosing $\bar{b}=1$:
\begin{equation}\label{besovthmpara}\aligned& 1< p\leq \tilde{p}<\infty
\\&1\leq q\leq \infty
\\&s> 1,~~ \bar{s}\tilde{p}<n
\\&0<2/\sigma=s-n/\tilde{p}-b<1
\\&0\leq \bar{s}< s-1
\\&1<\frac{n\tilde{p}}{2n-\bar{s}\tilde{p}}<\infty
\\& r>n/p
\\& s-1\leq \frac{n}{\tilde{p}}+r-n/p\leq s
\\& \sigma/2\leq p\leq \sigma.
\endaligned
\end{equation}
As in Section $4$ of \cite{sobpaper}, we can also obtain Theorem $\ref{special case 2}$ by modifying the standard method to fit the specific cases.  We also have the following result which extends the range of parameters allowed on the auxiliary space.
\begin{proposition}\label{higher regularity theorem2}Let $k>s_2>s_1$, with $s_2\geq 1$, and let $\varepsilon$ be a small positive number.  Then, for $k-s_2=s_2-s_1=\varepsilon$, for any solution $u$ to the LANS equation $(\ref{LANS})$ where 
\begin{equation*}u\in BC([0,T):B^{s_1}_{p,q}(\mathbb{R}^n)\cap L^{2/(s_2-s_1)}((0,T):B^{s_2}_{p,q}(\mathbb{R}^n)),
\end{equation*}
we have that $u\in L^1((0,T):B^{k}_{p,q}(\mathbb{R}^n))$ provided $s_2\geq n/p$.
\end{proposition}
We remark that this result can be extended, with the additional requirement that, for $k$ such that $k-s_1\geq 2$, $u\in L^a((0,T):B^{k}_{\tilde{p},q}(\mathbb{R}^n)$ requires $\tilde{p}>p$.

\bibliographystyle{amsplain}
\bibliography{references2}

\providecommand{\bysame}{\leavevmode\hbox to3em{\hrulefill}\thinspace}
\providecommand{\MR}{\relax\ifhmode\unskip\space\fi MR }
% \MRhref is called by the amsart/book/proc definition of \MR.
\providecommand{\MRhref}[2]{%
  \href{http://www.ams.org/mathscinet-getitem?mr=#1}{#2}
}
\providecommand{\href}[2]{#2}
\begin{thebibliography}{10}

\bibitem{chae}
D.~Chae, \emph{Local existence and blow-up criterion for the {E}uler equations
  in the {B}esov spaces}, Asymptotic {A}nalysis \textbf{38} (2004), 339--358.

\bibitem{CHMZ}
S.~Chen, D.~D. Holm, L.~Margolin, and R.~Zhang, \emph{Direct numerical
  simulations of the {N}avier-{S}tokes alpha model}, Phys. D \textbf{133}
  (1999), no.~1-4, 66--83, Predictability: quantifying uncertainty in models of
  complex phenomena (Los Alamos, NM, 1998).

\bibitem{MRS}
J.~Marsden, T.~Ratiu, and S.~Shkoller, \emph{The geometry and analysis of the
  averaged {E}uler equations and a new diffeomorphism group}, Geom. Funct.
  Anal. \textbf{10} (2000), no.~3, 582--599.

\bibitem{MS}
J.~Marsden and S.~Shkoller, \emph{Global well-posedness for the {L}agrangian
  averaged {N}avier-{S}tokes equations on bounded domains}, {P}hil. {T}rans.
  {R}. {S}oc. {L}ond. (2001), no.~359, 1449--1468.

\bibitem{MS2}
\bysame, \emph{The {A}nisotropic {L}agrangian {A}veraged {E}uler and
  {N}avier-{S}tokes equations}, {A}rch. {R}ational {M}ech. {A}nal. (2003),
  no.~166, 27--46.

\bibitem{MKSM}
K.~Mohseni, B.~Kosovi{\'c}, S.~Shkoller, and J.~Marsden, \emph{Numerical
  simulations of the {L}agrangian averaged {N}avier-{S}tokes equations for
  homogeneous isotropic turbulence}, Phys. Fluids \textbf{15} (2003), no.~2,
  524--544.

\bibitem{sobpaper}
N.~Pennington, \emph{Lagrangian {A}veraged {N}avier-{S}tokes equations with
  rough data in {S}obolev space}, http://arxiv.org/abs/1011.1856v2.

\bibitem{Shkoller}
S.~Shkoller, \emph{On incompressible averaged {L}agrangian hydrodynamics},
  E-print, (1999), http://xyz.lanl.gov/abs/math.AP/9908109.

\bibitem{SK}
\bysame, \emph{Analysis on groups of diffeomorphisms of manifolds with boundary
  and the averaged motion of a fluid}, Journal of {D}ifferential {G}eometry
  (2000), no.~55, 145--191.

\bibitem{taonde}
T.~Tao, \emph{Nonlinear {D}ispersive {E}quations}, American Mathematical
  Society, 2006.

\bibitem{TT}
M.~Taylor, \emph{Tools for {P}{D}{E}}, {M}athematical {S}urveys and
  {M}onographs, {V}ol. 81, {A}merican {M}athematical {S}ociety, Providence
  {R}{I}, 2000.

\end{thebibliography}

\end{document}